\documentclass[graybox]{svmult}

\usepackage{mathptmx}        
\usepackage{helvet}          
\usepackage{courier}         
\usepackage{type1cm}         
\usepackage{makeidx}         
\usepackage{graphicx}        
\usepackage{multicol}        
\usepackage[bottom]{footmisc}

\usepackage{amsmath}
\usepackage{amssymb}
\usepackage{amscd}
\usepackage{indentfirst}

\def\1{{\bf 1}}

\makeindex             

\usepackage[T1]{fontenc}                         
\usepackage{amsmath,amsfonts}                    
\DeclareMathAlphabet{\mathbbmsl}{U}{bbm}{m}{sl}  

\newtheorem{statement}{Statement}

\begin{document}

\title*{Double-sided Taylor's approximations and their applications in theory of trigonometric inequalities}
\titlerunning{Double-sided Taylor's approximations and their applications} 
\author{Branko Male\v sevi\' c${}^{\mbox{\scriptsize 1)}}$,
Tatjana Lutovac${}^{\mbox{\scriptsize 1)}}$,
Marija Ra\v sajski${}^{\mbox{\scriptsize 1)}}$,
Bojan Banjac${}^{\mbox{\scriptsize 2)}}$}
\authorrunning{B.~Male\v sevi\' c, T.~Lutovac, M.~Ra\v sajski, B. Banjac}
\institute{Branko Male\v sevi\' c, Tatjana Lutovac, Marija Ra\v sajski and Bojan Banjac
\at ${}^{1)}$School of Electrical Engineering, University of Belgrade; ${}^{2)}$Faculty of Technical Sciences, University of Novi Sad\\
Corresponding author \email{branko.malesevic@etf.bg.ac.rs}}

\maketitle

\abstract{In this paper  the double-sided {\sc Taylor}'s approximations are used to obtain
generalisations  and improvements of some trigonometric inequalities.}

\section{Introduction}

Many mathematical and engineering problems cannot be solved without {\sc Taylor}'s approximations
\cite{D_S_Mitrinovic_1970},
\cite{Milovanovic_Rassias_2014},
\cite{M_J_Cloud_B_C_Drachman_L_P_Lebedev_2014}.
Particularly, their application in proving various analytic inequalities is of great importance
\cite{C_Mortici_2011},
\cite{B_Malesevic_M_Makragic_JMI_2016},
\cite{Milica_Makragic_JMI_2017},
\cite{T_Lutovac_B_Malesevic_C_Mortici_JIA_2017}.
Recently, numerous inequalities have been generalized and improved by the use of the so-called
double-sided {\sc Taylor}'s approximations
\cite{Milica_Makragic_JMI_2017},
\cite{H_Alzer_M_K_Kwong_2017},
\cite{B_Malesevic_T_Lutovac_M_Rasajski_C_Mortici_Adv._Difference_Equ._2018}-\cite{M_Nenenzic_L_Zhu_AADM_2018}
and \cite{B_Malesevic_M_Rasajski_T_Lutovac_2019}.
Many topics regarding these approximations  are presented in
\cite{B_Malesevic_M_Rasajski_T_Lutovac_2019}.
Some of the basic concepts and results  about  the double-sided {\sc Taylor}'s approximations presented in \cite{B_Malesevic_M_Rasajski_T_Lutovac_2019}, which  will be used in this paper, are given in the
next section.

In this paper, using the double-sided {\sc Taylor}'s approximations, we obtain  generalizations and improvements of some trigonometric inequalities proved by {\sc J. Sandor} \cite{J_Sandor_2016}.

\smallskip

\begin{statement}
$(\mbox{\rm \cite{J_Sandor_2016}},~\mbox{\rm Theorem}~1)$
\label{J_Sandor_Th_01}
\begin{equation}
\label{Eq_01}
\displaystyle\frac{3}{8}
<
\displaystyle\frac{1-\mbox{\small $\displaystyle\frac{\cos x}{\cos \frac{x}{2}}$}}{x^2}
<
\displaystyle\frac{4}{\pi^2},
\end{equation}
for any $x \!\in\! (0,\pi/2)$.
\end{statement}

Note that {\sc J. D'Aurizio} \cite{J_D_Aurizio_2014} used the infinite products as well as some inequalities connected with
the {\sc Riemann} zeta function $\zeta$ to prove the right-hand side inequality (\ref{Eq_01}).

\medskip

\begin{statement}$(\mbox{\rm \cite{J_Sandor_2016}},~\mbox{\rm Theorem}~2)$
\label{J_Sandor_Th_02}
\begin{equation}
\label{Eq_02}
\displaystyle\frac{4}{\pi^2}\left(2 - \sqrt{2}\right)
<
\displaystyle\frac{2-\mbox{\small $\displaystyle\frac{\sin x}{\sin \frac{x}{2}}$}}{x^2}
<
\displaystyle\frac{1}{4},
\end{equation}
for any $x \!\in\! (0,\pi/2)$.
\end{statement}

Inequalities (\ref{Eq_01}) and (\ref{Eq_02}) are reducible to mixed trigonometric-polynomial inequalities
and can be proved by methods and algorithms that have been developed and shown in papers
\cite{B_Malesevic_M_Makragic_JMI_2016}, \cite{T_Lutovac_B_Malesevic_C_Mortici_JIA_2017}
and dissertation \cite{B_D_Banjac_2019}.

In this paper, we propose and prove generalizations of inequality (\ref{Eq_01}) by determining the sequence of the polynomial approximations. Also, an improvement of inequality (\ref{Eq_02}) is given for some intervals.
The proposed generalizations and improvements are based on the double-sided {\sc Ta\-ylo\-r}'s approximations and the corresponding results presented in \cite{B_Malesevic_M_Rasajski_T_Lutovac_2019}.

\section{An overview of the results related to double-sided \textsc{Ta\-ylo\-r}'s \break approximations}

Let us consider a real function $f : (a, b) \longrightarrow \mathbb{R}$,  such that there exist finite limits
 $f^{(k)}(a+)=\!\lim\limits_{x \rightarrow a+}{f^{(k)}(x)}$, for $k=0,1,\ldots,n$. \\
{\sc Taylor}'s polynomial
$$
T_n^{f,\,a+}(x)
=\displaystyle\sum_{k=0}^{n}{\displaystyle\frac{f^{(k)}(a+)}{k!}(x-a)^k}, \; n \!\in\! \mathbb{N}_0,
$$
and the polynomial
$$
\mbox{$\mathbbmsl{T}$}_n^{f;\,a+,\,b-}(x)
=
\left\{
\begin{array}{ccl}
T_{n-1}^{f,\,a+}(x)
+
\displaystyle\frac{1}{(b - a)^n}R_{n}^{f,\,a+}(b-)(x-a)^n                       &,& n \geq 1 \\[2.5 ex]
f(b-)                                                                           &,& n = 0,
\end{array}
\right.
$$
are called  the {\em first {\sc Taylor}'s approximation for the function $f$ in the right neighborhood~of~$a$},
and the {\em second {\sc Taylor}'s approximation for the function $f$ in the right neighborhood of $a$},
respectively.

Also,  the following functions:
$$
R_{n}^{f,\,a+}(x)
=
f(x) - T_{n-1}^{f,\,a+}(x), ~~ n \!\in\! \mathbb{N}
$$
and
$$
{\mbox{$\mathbbmsl{R}$}}_{n}^{f;\,a+,\,b-}(x)
=
f(x) - \mbox{$\mathbbmsl{T}$}_{n-1}^{f;\,a+,\,b-}(x), ~~ n \!\in\! \mathbb{N}
$$
are called the {\em remainder of the first {\sc Taylor}'s approximation in the right neighborhood of $a$},
and the {\em remainder of the second {\sc Taylor}'s approximation
in the right neighborhood of $a$}, respectively.

\medskip
The following Theorem,  which has been proved in \cite{S_Wu_L_Debnath_2009} and whose variants are considered
in \cite{S_Wu_HM_Srivastva_2008a}, \cite{S_Wu_L_Debnath_2008} and \cite{S_Wu_HM_Srivastva_2008b},  provides an
important result regarding {\sc Taylor}'s approximations.



\begin{theorem} \label{Theorem_1}
$(\mbox{\rm \cite{S_Wu_L_Debnath_2009}},~\mbox{\rm Theorem}~2)$
Suppose that $f(x)$ is a real function on $(a,b)$, and that $n$ is a positive integer such that $f ^{(k)}(a+)$,
for $k \!\in\! \{0,1,2, \ldots ,n\}$, exist.

\medskip
\noindent Supposing that $f^{(n)}(x)$ is increasing on $(a,b)$, then for all $x \!\in\! (a,b)$
the following inequality also holds$\,:$
\begin{equation}
\label{(2)}
 T_n^{f,\,a+}(x) < f(x) < \mbox{$\mathbbmsl{T}$}_n^{f;\,a+,\,b-}(x).
 \end{equation}
Furthermore, if $f^{(n)}(x)$ is decreasing on $(a,b)$, then the reversed inequality~of~\mbox{\rm (\ref{(2)})} holds.
\end{theorem}

The above theorem is called  {\em Theorem on double-sided \textsc{Taylor}'s approximations} in
\cite{B_Malesevic_M_Rasajski_T_Lutovac_2019}, i.e. {\em Theorem WD} in
\cite{B_Malesevic_T_Lutovac_M_Rasajski_C_Mortici_Adv._Difference_Equ._2018}-\cite{M_Nenenzic_L_Zhu_AADM_2018}.

\medskip
The proof of the following proposition is given in   \cite{B_Malesevic_M_Rasajski_T_Lutovac_2019}.

\begin{proposition} \label{Proposition_1}
$(\mbox{\rm \cite{B_Malesevic_M_Rasajski_T_Lutovac_2019}},~\mbox{\rm Proposition}~1)$
Consider a real function $f : (a, b) \longrightarrow \mathbb{R}$ such that there exist its first and second {\sc Taylor}'s approximations, for some $n \!\in\! {N}_0$. Then,
$$
\label{sgn}
{\mathop{\rm sgn}} {\Big (} \mbox{$\mathbbmsl{T}$}_{n}^{f,\, a+, \, b-}(x)\,-\,\mbox{$\mathbbmsl{T}$}_{n+1}^{f,\, a+, \, b-}(x) {\Big )}
=
{\mathop{\rm sgn}} {\Big (} {f(b-)\,-\,T_{n}^{f,\,a+}(b)} {\Big )},
$$
for all $ x\in(a,b)$.
\end{proposition}

From the above proposition, as shown in \cite{B_Malesevic_M_Rasajski_T_Lutovac_2019}, the following theorem
directly follows:

\begin{theorem} \label{Theorem_3}
$(\mbox{\rm \cite{B_Malesevic_M_Rasajski_T_Lutovac_2019}},~\mbox{\rm Theorem}~4)$
Consider the real analytic functions  $f : (a, b) \longrightarrow \mathbb{R}$:
$$
\label{f(x)}
f(x) = \sum_{k=0}^{\infty}{c_k(x-a)^k},
$$
where $c_k \!\in\! \mathbb{R}$ and  $c_k \geq 0$ for all $k \!\in\! \mathbb{N}_0$. Then,
 $$
 \begin{array}{c}
 T_0^{f,\,a+}(x) \leq \ldots\leq T_n^{f,\,a+}(x) \leq T_{n+1}^{f,\,a+}(x) \leq \ldots                     \\[1.0 ex]
 \ldots \leq  f(x) \leq \ldots                                                         \\[1.0 ex]
 \ldots \leq  \mbox{$\mathbbmsl{T}$}_{n+1}^{f;\,a+,\,b-}(x)
 \leq   \mbox{$\mathbbmsl{T}$}_n^{f;\,a+,\,b-}(x)
 \leq  \ldots  \leq \mbox{$\mathbbmsl{T}$}_0^{f;\,a+,\,b-}(x),
 \end{array}
 $$
for all $x \!\in\! (a,b)$.
\end{theorem}

\section{Main results}

$\,$ \vspace*{-3.5 mm}

\subsection{Generalization of Statement 1}

\smallskip

Consider the function:
$$
f(x)
=
\left\{
\begin{array}{ccc}
\mbox{\small $\displaystyle\frac{3}{8}$}                                                &,& x = 0,            \\[3.0 ex]
\displaystyle\frac{1-\mbox{\small $\displaystyle\frac{\cos x}{\cos \frac{x}{2}}$}}{x^2} &,& x \!\in\! (0,\pi).
\end{array}
\right.
$$

First, we prove that $ f $ is a real analytic function on  $ [0,\pi)$.
 Based on the elementary equality$:$
$$
1-\displaystyle\frac{\cos x}{\cos \frac{x}{2}}
=
1+\sec \frac{x}{2} - 2 \cos \frac{x}{2},
$$
and well known  power series expansions  \cite{I_Gradshteyn_I_Ryzhik_2014} (formula 1.411):

$$
\begin{array}{ccc}
\cos t = \displaystyle\sum_{k=0}^{\infty}{\mbox{\small $\displaystyle\frac{(-1)^k}{(2k)!}$} \,t^{2k}}  & \qquad  & t \!\in\! R,      \\[5.0 ex]
\sec t = \displaystyle\sum_{k=0}^{\infty}{\mbox{\small $\displaystyle\frac{|E_{2k}|}{(2k)!}$}\,t^{2k}} & \qquad  &
t \!\in\! \left(-\mbox{\small $\displaystyle\frac{\pi}{2}$},\mbox{\small $\displaystyle\frac{\pi}{2}$}\right);
\end{array}
$$

\noindent
where $E_k$ are {\sc Euler}'s numbers \cite{I_Gradshteyn_I_Ryzhik_2014},  for
$t=\mbox{\small $\displaystyle\frac{x}{2}$} \!\in\! \left[0,\mbox{\small $\displaystyle\frac{\pi}{2}$}\right)$,
i.e. for $x \!\in\! [0,\pi)$, we have$:$
$$
f(x) = \displaystyle\sum_{k=1}^{\infty}{\mbox{\small $\displaystyle\frac{|E_{2k}| - 2(-1)^k}{2^{2k}(2k)!}$}x^{2k-2}}
$$

\vspace*{-2.0 mm}

\noindent
i.e.

\vspace*{-1.5 mm}

$$
f(x)
=
\mbox{\small $\displaystyle\frac{3}{8}$}
+
\mbox{\small $\displaystyle\frac{1}{128}$}\,x^2
+
\mbox{\small $\displaystyle\frac{7}{5120}$}\,x^4
+
\mbox{\small $\displaystyle\frac{461}{3440640}$}\,x^6
+
\mbox{\small $\displaystyle\frac{16841}{1238630400}$}\,x^8
+
\ldots
$$

\smallskip\noindent
where the power series converges for $x \!\in\! [0,\pi)$.

\smallskip
Further, based on the elementary well-known features of
 {\sc Euler}'s numbers $E_k$, we have$:$
$$
c_{2k-2} = \mbox{\small $\displaystyle\displaystyle\frac{|E_{2k}| - 2(-1)^k}{2^{2k}(2k)!}$} > 0
\quad\mbox{and}\quad
c_{2k-1} = 0,
$$
for $k = 1, 2, ...\,$.

\smallskip
Finally, from Theorem~\ref{Theorem_3} the following result directly follows.

\break

\begin{theorem} \label{Generalization_Statement_1}
For the function
$$
f(x)
=
\left\{
\begin{array}{ccc}
\mbox{\small $\displaystyle\frac{3}{8}$}                                                &,& x = 0,            \\[2.0 ex]
\displaystyle\frac{1-\mbox{\small $\displaystyle\frac{\cos x}{\cos \frac{x}{2}}$}}{x^2} &,& x \!\in\! (0,\pi)
\end{array}
\right.
$$
 and any  $c\in (0, \pi)$  the following inequalities hold true:
\begin{equation}
\label{generalizacija_Shandor_1}
\begin{array}{c}
\mbox{\small $\displaystyle\frac{3}{8}$}
= T_0^{f,\,0+}(x) \leq T_2^{f,\,0+}(x)\leq \ldots\leq T_{2n}^{f,\,0+}(x) \leq \ldots                      \\[2.0 ex]
 \ldots \leq  f(x) \leq \ldots                                                                            \\[1.0 ex]
 \leq   \mbox{$\mathbbmsl{T}$}_{2m}^{f;\,0+,\,c-}(x)
 \leq  \ldots \leq \mbox{$\mathbbmsl{T}$}_2^{f;\,0+,\,c-}(x) \leq \mbox{$\mathbbmsl{T}$}_0^{f;\,0+,\,c-}(x)
= \left(1-\mbox{\small $\displaystyle\frac{\cos c}{\cos \frac{c}{2}}$}\right)\mbox{\Large $/$}{c^2}.
\end{array}
\end{equation}
for every $x \!\in\! \left(0, c\right)$, where  $m, n \!\in\! \mathbb{N}_0$.

\end{theorem}

\medskip
\noindent Note that inequalities from Statement  1 can be  directly obtained from (\ref{generalizacija_Shandor_1}),
for~$c\!=\!\mbox{\small $\displaystyle\frac{\pi}{2}$}$

$$
\mbox{\small $\displaystyle\frac{3}{8}$}
<
\displaystyle\frac{1-\mbox{\small $\displaystyle\frac{\cos x}{\cos \frac{x}{2}}$}}{x^2}
<
\mbox{\small $\displaystyle\frac{4}{\pi^2}$}.
$$

Also,  Theorem~\ref{Generalization_Statement_1} gives a generalization and a sequence of improvements of
results from Statement 1. For example, for $c=\pi/2$ i.e. for $x \!\in\! \left(0, \mbox{\small $\displaystyle\frac{\pi}{2}$}\right)$
we have:
$$
\mbox{\small $\displaystyle\frac{3}{8}$}
\leq
T_{2}^{f,\,0+}\!(x)
\!=\!
\mbox{\small $\displaystyle\frac{3}{8}$}
\!+\!
\mbox{\small $\displaystyle\frac{1}{128}$}x^2
\leq
f(x)
\leq
\mbox{$\mathbbmsl{T}$}_{2}^{f;\,0+,\,\pi/2-}\!(x)
\!=\!
\mbox{\small $\displaystyle\frac{3}{8}$}
\!+\!
\left(\mbox{\small $\displaystyle\frac{16}{\pi^4}$} \!-\! \mbox{\small $\displaystyle\frac{3}{2 \pi^2}$}\right)x^2
\leq
\mbox{\small $\displaystyle\frac{4}{\pi^2}$}.
$$
Using standard numerical methods it is easy to verify:
$$
\mathop{max}\limits_{x \in [0,\pi/2]}{|R_3^{f,0+}(x)|}
=
f(\pi/2)
=
0.01100 ...
$$
and
$$
\mathop{max}\limits_{x \in [0,\pi/2]}|\mbox{$\mathbbmsl{R}$}_{3}^{f;\,0+,\,\pi/2-}(x)|
=
f(1.14909...)
=
0.00315 ... \, .
$$

\vspace*{-2.0 mm}

\subsection{An improvement of  Statement 2}

Let $\beta \!\in\! (0,\pi)$ be a fixed real number.
Consider the function$:$
$$
g(x)
=
\left\{
\begin{array}{ccc}
\mbox{\small $\displaystyle\frac{1}{4}$}                                                &,& x = 0,            \\[2.0 ex]
\displaystyle\frac{1-\mbox{\small $\displaystyle\frac{\sin x}{\sin \frac{x}{2}}$}}{x^2} &,& x \!\in\! (0,\beta].
\end{array}
\right.
$$

\break

We prove that $g$ is a real analytic function on  $[0,\beta]$.

\smallskip
\noindent
Notice  that
\begin{equation}
\label{g1-g2}
g(x) = g_1(x) - g_2(x)
\end{equation}
for $x \!\in\! [0,\beta]$, where
$$
g_1(x)
=
\left\{
\begin{array}{ccc}
\mbox{\small $\displaystyle\frac{1}{4}$}                            &,& x = 0,            \\[2.25 ex]
\displaystyle\frac{\cosh\!\frac{x}{2} - \cos\!\frac{x}{2}}{x^2}     &,& x \!\in\! (0,\beta]
\end{array}
\right.
$$
and
$$
g_2(x)
=
\left\{
\begin{array}{ccc}
0                                                                    &,& x = 0,            \\[2.25 ex]
\displaystyle\frac{ \cosh\!\frac{x}{2} + \cos\!\frac{x}{2} - 2}{x^2} &,& x \!\in\! (0,\beta].
\end{array}
\right.
$$

Since the functions $g_1$ and $g_2$ are real analytic functions  on $[0,\beta]$, with the following power series expansions:
$$
g_1(x)
=
\displaystyle\sum_{k=0}^{\infty}{\mbox{\small $\displaystyle\frac{1}{2^{4k}(4k+2)!}$}x^{4k}}
$$
and
$$
g_2(x)
=
\displaystyle\sum_{k=0}^{\infty}{\mbox{\small $\displaystyle\frac{1}{2^{4k+2}(4k+4)!}$}x^{4k+2}},
$$
the function $g$   must also be a real analytic function on $[0,\beta]$.

\medskip
Also, from Theorem~\ref{Theorem_3} the following results directly follow.
\begin{theorem}
\label{g1}
For all $c\in \left(0, \pi\right)$  the following inequalities hold true$:$
$$
\begin{array}{c}
\mbox{\small $\displaystyle\frac{1}{4}$}
= T_0^{g_1,\,0+}(x) \leq \ldots\leq T_{4n}^{g_1,\,0+}(x) \leq T_{4n+4}^{g_1,\,0+}(x) \leq \ldots        \\[1.50 ex]
 \ldots \leq  g_1(x) \leq \ldots                                                                        \\[1.25 ex]
 \ldots \leq  \mbox{$\mathbbmsl{T}$}_{4m+4}^{g_1;\,0+,\,c-}(x)
 \leq   \mbox{$\mathbbmsl{T}$}_{4m}^{g_1;\,0+,\,c-}(x)
 \leq  \ldots  \leq \mbox{$\mathbbmsl{T}$}_0^{g_1;\,0+,\,c-}(x)
= g_1(c).
\end{array}
$$
for all $x \!\in\! \left(0, c\right)$,   where $m, n \!\in\! \mathbb{N}_0$.
\end{theorem}
\begin{theorem}
\label{g2}
For all $c\in \left(0, \pi\right)$  the following inequalities hold true$:$
$$
\begin{array}{c}
\mbox{\small $\displaystyle\frac{1}{192}$}x^2
=T_2^{g_2,\,0+}(x)
\leq
\ldots
\leq
T_{4n+2}^{g_2,\,0+}(x)
\leq
T_{4n+6}^{g_2,\,0+}(x)  \leq \ldots                    \\[1.50 ex]
\ldots \leq g_2(x) \leq \ldots                         \\[1.00 ex]
\ldots \leq
\mbox{$\mathbbmsl{T}$}_{4m+6}^{g_2;\,0+,\,c-}(x)
\leq
\mbox{$\mathbbmsl{T}$}_{4m+2}^{g_2;\,0+,\,c-}(x)
\leq \ldots
\leq
\mbox{$\mathbbmsl{T}$}_2^{g_2;\,0+,\,c-}(x)
=
\mbox{\small $\displaystyle\frac{g_2(c)}{c^2}x^2$}
\end{array}
$$
for all $x \!\in\! \left(0, c \right)$, where $m, n \!\in\! \mathbb{N}_0$.
\end{theorem}

Thus, from  (\ref{g1-g2}), Theorem~\ref{g1} and Theorem~\ref{g2},  for $c=\mbox{\small $\displaystyle\frac{\pi}{2}$}$, an improvement of inequalities from Statement  2 are obtained, as shown bellow.

\medskip
First, for all $x \!\in\! \left(0, \pi/2\right)$  the following inequalities hold true$:$

$$
\mbox{\small $\displaystyle\frac{1}{4}$}
-
\mbox{\small $\displaystyle\frac{4}{\pi^2}$}\,g_2\!\!\left(\mbox{\small $\displaystyle\frac{\pi}{2}$}\right)x^2
\leq
g(x)
\leq
g_1\!\!\left(\mbox{\small $\displaystyle\frac{\pi}{2}$}\right) - \mbox{\small $\displaystyle\frac{1}{192}$} \, x^2
$$
i.e.
$$
\mbox{\small $\displaystyle\frac{1}{4}$}
-
\mbox{\small $\displaystyle\frac{16}{\pi^4}$}\displaystyle\left( \cosh\!\mbox{\small $\displaystyle\frac{\pi}{4}$}
+ \mbox{\small $\displaystyle\frac{\sqrt{2}}{2}$} - 2\right)x^2
\leq
g(x)
\leq
\displaystyle \mbox{\small $\displaystyle\frac{4}{\pi^2}$}\!\left(
\cosh\!\mbox{\small $\displaystyle\frac{\pi}{4}$} - \mbox{\small $\displaystyle\frac{\sqrt{2}}{2}$}\right)
-
\mbox{\small $\displaystyle\frac{1}{192}$} \, x^2.
$$

It is easy to check
$$
\mbox{\small $\displaystyle\frac{4}{\pi^2}$}\!\left(
\cosh\!\mbox{\small $\displaystyle\frac{\pi}{4}$}
-
\mbox{\small $\displaystyle\frac{\sqrt{2}}{2}$}\right)
-
\mbox{\small $\displaystyle\frac{1}{192}$} \, x^2
\leq
\mbox{\small $\displaystyle\frac{1}{4}$}
$$
for all $x \!\in\! [\delta_2, \pi/2]$, where
$\delta_2 = \mbox{\small $\displaystyle\frac{4 \sqrt{3} e^{-\frac{\pi}{8}}}{\pi}$}
\sqrt{8 + 8 e^{\frac{\pi}{2}} - (\pi^2+8\sqrt{2})e^{\frac{\pi}{4}}} = 0.22525 ...\,$.

\smallskip Also,
$$
\mbox{\small $\displaystyle\frac{4}{\pi^2}$}(2-\sqrt{2})
\leq
\mbox{\small $\displaystyle\frac{1}{4}$}
-
\mbox{\small $\displaystyle\frac{16}{\pi^4}$}\!\left(
\cosh\!\mbox{\small $\displaystyle\frac{\pi}{4}$}
+
\mbox{\small $\displaystyle\frac{\sqrt{2}}{2}$}
-
2\right)x^2
$$
for all $x \!\in\! [0,\delta_1]$, where $\delta_1
=
\mbox{\small $\displaystyle\frac{\sqrt{2} \pi e^{\frac{\pi}{8}}\sqrt{\pi^2+16\sqrt{2}-32}}{8\sqrt{(\sqrt{2}-4)e^{\frac{\pi}{4}}
+e^{\frac{\pi}{2}}+1}}$} =1.55456 ... \,$.

\section{Conclusion}

In this paper, we showed a way to prove some trigonometric inequalities using the double-sided {\sc Taylor}'s approximations.
The presented approach enabled generalizations of inequalities (\ref{Eq_01}) i.e. produced sequences of polynomial approximations
of the given trigonometric function $f$.

Note that  Theorem~\ref{Theorem_3} cannot be applied directly to inequality (\ref{Eq_02}) because the function $g$ has
an alternating series expansion. We overcame this obstacle by representing this function by a linear combination of
two functions whose power series expansions have nonnegative coefficients.

Our approach  makes a good basis for the systematic proving of trigonometric inequalities. Developing general,
automated-oriented methods for proving of trigonometric inequalities is an area our continuing interest
\cite{B_Banjac_M_Nenenzic_B_Malesevic_Telfor_2015},
\cite{B_Malesevic_M_Makragic_JMI_2016}-\cite{B_Banjac_M_Makragic_B_Malesevic_Results_2016},
\cite{Milica_Makragic_JMI_2017}-\cite{B_Malesevic_I_Jovovic_B_Banjac_JMI_2017},
\cite{B_Malesevic_T_Lutovac_M_Rasajski_C_Mortici_Adv._Difference_Equ._2018}-\cite{B_Malesevic_M_Rasajski_T_Lutovac_2019}
and \cite{B_D_Banjac_2019}.

\bigskip
\textbf{Acknowledgment.}
Research of the first and second and third author was supported in part by the Serbian Ministry of
Education, Science and Technological Development, under Projects ON 174032 \& III 44006, ON 174033
and TR 32023, respectively.

\end{document}